\documentclass[reqno,11pt]{amsart}
\usepackage{amscd,amssymb}

\theoremstyle{plain}
\newtheorem{Thm}[subsection]{Theorem}
\newtheorem{Cor}[subsection]{Corollary}
\newtheorem{Lem}[subsection]{Lemma}
\newtheorem{Prop}[subsection]{Proposition}
\newtheorem{Conj}[subsection]{Conjecture}

\theoremstyle{definition}
\newtheorem{Def}[subsection]{Definition}

\theoremstyle{remark}

\newtheorem{Rem}[subsection]{Remark}

\newtheorem{conj}{Conjecture}

\errorcontextlines=0 \numberwithin{equation}{section}
\renewcommand{\rm}{\normalshape}

\newif\ifShowLabels
\ShowLabelstrue
\newdimen\theight
\def\TeXref#1{%
    \leavevmode\vadjust{\setbox0=\hbox{{\tt
        \quad\quad  {\small \rm #1}}}%
    \theight=\ht0
    \advance\theight by \lineskip
    \kern -\theight \vbox to
    \theight{\rightline{\rlap{\box0}}%
    \vss}%
    }}%

\ShowLabelsfalse

\renewcommand{\sec}[2]{\section{#2}\label{S:#1}%
    \ifShowLabels \TeXref{{S:#1}} \fi}
\newcommand{\ssec}[2]{\subsection{#2}\label{SS:#1}%
    \ifShowLabels \TeXref{{SS:#1}} \fi}

\newcommand{\refs}[1]{Section ~\ref{S:#1}}
\newcommand{\refss}[1]{Section ~\ref{SS:#1}}

\newcommand{\reft}[1]{Theorem ~\ref{T:#1}}
\newcommand{\refl}[1]{Lemma ~\ref{L:#1}}
\newcommand{\refp}[1]{Proposition ~\ref{P:#1}}

\newcommand{\refe}[1]{\eqref{E:#1}}
\newcommand{\refco}[1]{Conjecture ~\ref{Co:#1}}

\newenvironment{thm}[1]%
    { \begin{Thm} \label{T:#1}  \ifShowLabels \TeXref{T:#1} \fi }%
    { \end{Thm} }

\renewcommand{\th}[1]{\begin{thm}{#1} \sl }
\renewcommand{\eth}{\end{thm} }

\newenvironment{lemma}[1]%
    { \begin{Lem} \label{L:#1}  \ifShowLabels \TeXref{L:#1} \fi }%
    { \end{Lem} }
\newcommand{\lem}[1]{\begin{lemma}{#1} \sl}
\newcommand{\elem}{\end{lemma}}

\newenvironment{propos}[1]%
    { \begin{Prop} \label{P:#1}  \ifShowLabels \TeXref{P:#1} \fi }%
    { \end{Prop} }
\newcommand{\prop}[1]{\begin{propos}{#1}\sl }
\newcommand{\eprop}{\end{propos}}

\newenvironment{corol}[1]%
    { \begin{Cor} \label{C:#1}  \ifShowLabels \TeXref{C:#1} \fi }%
    { \end{Cor} }
\newcommand{\cor}[1]{\begin{corol}{#1} \sl }
\newcommand{\ecor}{\end{corol}}

\newenvironment{defeni}[1]%
    { \begin{Def} \label{D:#1}  \ifShowLabels \TeXref{D:#1} \fi }%
    { \end{Def} }
\newcommand{\defe}[1]{\begin{defeni}{#1} \sl }
\newcommand{\edefe}{\end{defeni}}

\newenvironment{remark}[1]%
    { \begin{Rem} \label{R:#1}  \ifShowLabels \TeXref{R:#1} \fi }%
    { \end{Rem} }
\newcommand{\rem}[1]{\begin{remark}{#1}}
\newcommand{\erem}{\end{remark}}

\newenvironment{conjec}[1]%
    { \begin{Conj} \label{Co:#1}  \ifShowLabels \TeXref{Co:#1} \fi }%
    { \end{Conj} }
\renewcommand{\conj}[1]{\begin{conjec}{#1} \sl }
\newcommand{\econj}{\end{conjec}}

\newcommand{\eq}[1]%
    { \ifShowLabels \TeXref{E:#1} \fi
       \begin{equation} \label{E:#1} }
\newcommand{\eeq}{ \end{equation} }

\newcommand{\prf}{ \begin{proof} }
\newcommand{\epr}{ \end{proof} }

\newcommand\alp{\alpha}     

\newcommand\gam{\gamma}     \newcommand\Gam{\Gamma}
\newcommand\del{\delta}     \newcommand\Del{\Delta}

\newcommand\kap{\kappa}
\newcommand\lam{\lambda}        \newcommand\Lam{\Lambda}

\newcommand\ome{\omega}     \newcommand\Ome{\Omega}

\newcommand\calA{{\mathcal{A}}}
\newcommand\calB{{\mathcal{B}}}
\newcommand\calC{{\mathcal{C}}}
\newcommand\calD{{\mathcal{D}}}

\newcommand\calF{{\mathcal{F}}}
\newcommand\calG{{\mathcal{G}}}
\newcommand\calH{{\mathcal{H}}}
\newcommand\calI{{\mathcal{I}}}

\newcommand\calK{{\mathcal{K}}}
\newcommand\calL{{\mathcal{L}}}
\newcommand\calM{{\mathcal{M}}}

\newcommand\calO{{\mathcal{O}}}
\newcommand\calP{{\mathcal{P}}}

\newcommand\calR{{\mathcal{R}}}
\newcommand\calS{{\mathcal{S}}}
\newcommand\calT{{\mathcal{T}}}
\newcommand\calU{{\mathcal{U}}}

        \newcommand\bfE{{\mathbf E}}
        \newcommand\bfF{{\mathbf F}}
        \newcommand\bfG{{\mathbf G}}
        \newcommand\bfH{{\mathbf H}}
        \newcommand\bfI{{\mathbf I}}
        
        \newcommand\bfK{{\mathbf K}}

        \newcommand\bfS{{\mathbf S}}

\newcommand\TT{\mathbb{T}}

\renewcommand\SS{\mathbb{S}}

\newcommand\FF{\mathbb{F}}
\newcommand\GG{\mathbb{G}}
\newcommand\HH{\mathbb{H}}

\newcommand\ZZ{\mathbb{Z}}

\newcommand\CC{\mathbb{C}}
\newcommand\VV{\mathbb{V}}

 \newcommand\grg{{\mathfrak{g}}}

 \newcommand\grm{{\mathfrak{m}}}

\newcommand\sdp{\times \hskip -0.3em {\raise 0.3ex
\hbox{$\scriptscriptstyle |$}}} 

\newcommand\End{\operatorname{End\,}}

\newcommand\Gr{\operatorname{Gr}}

\newcommand\id{\operatorname{id}}

\newcommand\Ind{\operatorname{Ind}}

\newcommand\pr{\operatorname{pr}}

\newcommand\Spec{\operatorname{Spec}}

\newcommand\oG{{\overline{G}}}

\newcommand\hatG{{\widehat{G}}}

\newcommand\tilG{{\widetilde{G}}}

\newcommand\tils{{\widetilde{s}}}

\newcommand\tilY{{\widetilde{Y}}}

\newcommand\x{\times}
\newcommand\ten{\otimes}

\newcommand{\ra}{\rangle}
\newcommand{\la}{\langle}

\newcommand\Rep{\operatorname{Rep}}
\newcommand\hGG{\hat{\GG}}
\newcommand\Vect{\VV ect}
\newcommand\Inv{\operatorname{Inv}}
\renewcommand\Spec{\operatorname{Spec}}
\newcommand\hsp{\calH^{\text{sph}}}

\newcommand\Bun{\operatorname{Bun}}
\newcommand\tbg{\widetilde{\Bun_G}}
\newcommand\Pic{\operatorname{Pic}}

\newcommand\tgg{\widetilde{\GG}}
\newcommand\oGr{\overline \Gr}
\newcommand\tGr{\widetilde \Gr}
\newcommand\Fl{\calF l}
\newcommand\oFl{\overline \Fl}
\newcommand{\wbg}{{\widetilde \Bun_G}}
\newcommand\triv{\operatorname{triv}}
\newcommand\Pro{\operatorname{Pro}}
\newcommand\Coinv{\operatorname{Coinv}}
\newcommand\hspp{\calH^{\text{sph}'}}
\begin{document}
Dedicated to G.~Lusztig on the occasion of his 60th birthday

\bigskip

\title{Some examples of Hecke algebras for two-dimensional local fields}
\author{Alexander Braverman and David Kazhdan}
\address{Einstein Institute of Mathematics,
Edmond J.~Safra campus, Hebrew University of Jerusalem,
Givat-Ram, Jerusalem}
\begin{abstract}
Let $\bfK$ be a local non-archimedian field, $\bfF=\bfK((t))$
and let $\bfG$ be a split
semi-simple group.
The purpose of this paper is to study certain analogs of spherical and
Iwahori Hecke algebras for representations of the group $\GG=G(\bfF)$ and its
central extension $\hGG$ .
For instance our spherical Hecke algebra corresponds to the subgroup
$G(\calA)\subset G(\bfF)$ where $\calA\subset \bfF$ is the subring $\calO_{\bfK}((t))$
where $\calO_\bfK\subset \bfK$ is the ring of integers. It turns out that for generic
level (cf. \cite{GK1}) the  spherical Hecke algebra is trivial; however, on the critical level it is quite
large. On the other hand we expect that the
 size of the corresponding Iwahori-Hecke algebra does not depend on a choice of a level
 (details will be considered in another publication).
\end{abstract}
\maketitle
\sec{int}{Introduction}
\ssec{notation}{}Let Let $\bfK$ be a local non-archimedian field and let $\bfF=\bfK((t))$.
In this paper we shall actually assume that $\bfK=\FF_q((x))$ though this is
probably not necessary.
Denote by $\calO_{\bfK}$ the ring of integers of $\bfK$. Set also $\calO_{\bfF}=\bfK[[t]]$,
$\grm_{\bfF}=t\calO_{\bfF}$,
$\calA=\calO_{\bfK}((t))$, $\bfE=\FF_q((t))$ and $\calO_{\bfE}=\FF_q[[t]]$.
We denote by $Vect$ the category of vector spaces over $\CC$ and by $\Vect$ the category
of pro-vector spaces over $\CC$.
Let now $G$ be a connected split reductive algebraic group (defined over $\ZZ$). Set $\GG=G(\bfF)$.
The category of $\Rep(\GG)$ of representations of $\GG$ was defined in \cite{GK1}.
The group $\GG$ admits canonical central extension $\hGG$ by means of $\bfK^*$ and we denote by
$\Rep(\hGG)$ the category of its representations. For each $\kap:\bfK^*\to\CC^*$
we denote by $\Rep_{\kap}(\GG)$ the full subcategory of of $\Rep(\hGG)$ on which the central
$\bfK^*$ acts by $\kap$.
\ssec{}{Cherednik algebra as a Hecke algebra}In this paper we want to discuss some examples of ``Hecke algebras'' for the group
$\hGG$. One of such examples is discussed in detail in \cite{K2} and \cite{GK1}.
Namely, let us choose a Borel subgroup $B$ of $G$ defined over $\bfK$ and
let $U$ denote its unipotent radical and choose a Cartan subgroup $T$
in $B$. Let $\bfI\subset G(\calO_{\bfF}$ be the Iwahori subgroup consisting
of all elements of $G(\calO_{\bfF})$ whose reduction mod t  lies in $B(\bfK)$. Let also
$\bfI^0$ denote the preimage of $U(\bfK)$ under the map $\bfI\to B(\bfK)$ and set
$\bfI^{00}=T(\calO_{\bfF})$. Consider the functor
$$
\operatorname{Coinv}_{\bfI^{00}}:\Rep(\hGG)\to \Vect
$$
sending every representation to its coivariants with respect to $\bfI^{00}$.
Let $\overset{\cdot\cdot}\bfH$ denote the algebra of endomorphisms of this
functor. It is shown in \cite{GK2} that this algebra is naturally
isomorphic to Cherednik's double affine Hecke algebra associated with $G$.
This is an extension of a previous result of
Kapranov (cf. \cite{K2}).
\ssec{kolcoa}{The ring $\calA$ and the corresponding subgroups}
One of the purposes of the present paper is to try to understand
the results of \cite{K2} in some "global" context. In other words,
one would like to develop some kind of two-dimensional theory
of automorphic forms. Some speculations about this are presented
in the last section of this article; the main part of the paper
is devoted to a discussion of certain ``Hecke algebras'' which
are supposed to play the same role in this (not yet constructed) theory
as the usual Hecke algebras of p-adic groups play in the theory
of automorphic forms. In the theory of representations of the
group $G(\bfK)$ the Hecke algebras are attached to open compact
subgroups of this group; very often such subgroups can
be realized as subgroups of $G(\calO_{\bfK})$ of finite
index. Given such a group $\Gam$ we can consider the
functor $\Inv_{\Gam}:\Rep_{G(\bfK)}\to Vect$
sending every representation $V$ to the subspace $V^{\Gam}\subset V$
of $\Gam$-invariants. (Note that since $\Gam$ is compact we have
 $V^{\Gam}=V_{\Gam}$ where $V_{\Gam}$ is the quotient space of
 $\Gam$-coinvariants). The corresponding Hecke
algebra $H(G,\Gam)$ can be defined as the algebra of endomorphisms
of this functor.
It turns out (cf. some explanations in \refss{global}) that in our
case the relevant subgroups of $G(\bfF)$ are of the following
form. Let $\calA=\calO_{\bfK}((t))$.
\footnote{It is easy to see that $\calA$ (as a subring of $\bfF$)
does not depend on the choice of $t$}.
We have the natural homomorphism $\calA\to \FF_q((t))$ (reduction
modulo $\grm_{\bfK}$). Hence we get the natural homomorphism
$\eta:G(\calA)\to G(\FF_q((t)))=G(\bfE)$. The subgroups $\Gam\subset \GG$
in which we shall be interested in in this paper are those which
are equal to the inverse image of a closed subgroup of
$G(\bfE)$.
\footnote{One can consider a more general class of subgroups considering reductions
modulo higher powers of $\grm_{\bfK}$. However, we are not going to do it
in this paper.}
More specifically we are going to concentrate
on the following two examples:

1) $\Gam=G(\calA)$.

2) $\Gam=\bfI_{\calA}^{00}$ where
$$
\bfI_{\calA}^{00}=\eta^{-1}(T(\calO_{\bfE})U(\bfE)).
$$

We shall refer to the first case as the {\it spherical case} and to the
second as {\it the Iwahori case}.
It is easy to see that the central extension $\hGG$ splits over $G(\calA)$.
Since $G(\calA)$ is equal to it's commutator  we may consider the group $G(\calA)$  as
a subgroup of $\hGG$.
\ssec{}{The spherical Hecke algebra}.
Fix $\kap:\bfF^*\to \CC^*$ as above and assume that it is trivial on $\calO_{\bfF}^*$. In this case
we may write $\kap(x)=q^{c v(x)}$ where $c\in \CC$. In this case we shall write
$\Rep_c(\hGG)$ instead of $\Rep_\kap(\hGG)$.
Let $\Inv_c: \Rep_c(\hGG)\to \Vect$ be the functor of invariants with respect to $G(\calA)$
\footnote{Since we are dealing with pro-vector spaces the notion of invariants is a bit tricky (cf.
\refss{invc})}.
Set
$$
\hsp_c=\End(\Inv_c).
$$
In other words, $\hsp_c$ is the universal algebra acting on the (pro)space of $G(\calA)$-invariants in any
representation of level $c$.

\medskip
\noindent
{\bf Remarks.}

\noindent
1. In this case the functor of invariants no longer
coincides with the functor of coinvariants. It turns out that we
have to use the former (it is not difficult to see that the latter
is almost always trivial).

\smallskip
\noindent
2. One can also study another version of the "double spherical Hecke algebra" when the subgroup
$G(\calA)$ is replaced by $G(\calO_{\bfF})$ (and the functor of invariants is replaced
by the functor of coinvariants). This is somewhat simpler but seems to be less relevant for
"global" purposes discussed in \refss{global}.
\conj{sph-general}
\begin{enumerate}
\item
For any $c$ the algebra $\hsp_c$ is commutative.
\item
One has $\hsp_c=\CC$ unless $q^{c+h^{\vee}}$ is a root of unity; here $h^{\vee}$ denotes the dual Coxeter number of $G$.
\end{enumerate}
\econj
\noindent
{\bf Remark.} The reader should compare this Conjecture with the corresponding (known) statements
for affine Lie algebras. For example the analog of 1) for affine Lie algebras is a very general
statement which is proved using the theory of vertex operator algebras. It seems that in our case
one has to develop some sort of similar theory in order to prove such statements.

\bigskip
\noindent
\ssec{crit}{The case of critical level}
The size of the algebra $\hsp_c$ changes drastically when $q^{c+h^{\vee}}$ is a root of unity.
In this paper we shall actually restrict ourselves to the critical level, i.e. to the case when
 $c=-h^{\vee}$ (in this
case we shall write $\hsp_{crit}$ instead of $\hsp_{-h^{\vee}}$). While we still
can't show that in this case the algebra is commutative, we can construct a very large
commutative subalgebra in it. This is done as follows.
Let $\calB=\calO_{\bfF}\cap \calA=\FF_q[[x,t]]$.
Let us consider the scheme $\Del=\Spec(\calB)$. By an {\it irreducible curve}
in $\Del$ we shall mean a proper principal prime ideal of $\calB$. If such a curve
$C$ is given by the equation $f=0$ for some $f\in \calB$ (i.e. $f$ is a generator of the
corresponding ideal) then we say that $C$ is good if the ring $\calB$ is topologically
generated by $t$ and $f$ over $\FF_q$ (thus automatically we have
$\calB=\FF_q[[f,t]]$). In other words $C$ is a smooth irreducible curve
which intersects the curve $X:=\{ t=0\}$ transversely.
 Let $\calC$
denote the set of good curves. Let $\calR$ be the free commutative $\CC$-algebra
whose generators are elements of $\calC$. Let us also denote by $\Lam$ the coweight lattice
of $G$ and by $\Lam_+\subset \Lam$ the subsemigroup of dominant coweights.
\th{main}The algebra $\hsp_{crit}$ contains $\calR[\Lam_+]$ as a subalgebra.
\eth
In this section we give a rough  idea of the proof of \reft{main}
(which will also explain why the critical level is special) and provide details later in the paper.
The theorem says basically that given any $C\in \calC$ we may construct an embedding
$\alp_C:\CC[\Lam_+]\hookrightarrow \hsp_{crit}$ and that their images commute for different curves $C$.
In other words, given $C$ as above we need to construct
the corresponding element $\alp_C(\lam)\in \hsp_{crit}$ for each $\lam\in\Lam_+$ such that for any
$\lam,\mu\in\Lam_+$ we have
$$
\alp_C(\lam)\alp_C(\mu)=\alp_C(\lam+\mu)
$$
and such that for any two curves $C,C'$ and any $\lam,\mu\in\Lam_+$ the elements $\alp_C(\lam)$ and $\alp_{C'}(\mu)$
commute.

Let $C$ be given by the equation $f=0$. Choose a maximal torus $T$ in $G$ and
set $f^{\lam}=\lam(f)$.
Let $\VV$ be a representation of $\hGG$ of some level $c$. We need to construct
an operator $\alp_C^{\VV}:\VV^{G(\calA)}\to \VV^{G(\calA)}$.
Assume for simplicity that $\VV^{G(\calA)}$ is a vector space (and not
a pro-space) so that we can talk about its elements.
Then for
any $v\in \VV^{G(\calA)}$ we would like to write
\eq{alpC}
\alp_C^{\VV}(\lam)(v)= \int\limits_{G(\calA)f^{\lam}G(\calA)/G(\calA)} g(v) dg.
\end{equation}
For this formula to make sense we have to

\medskip
\noindent
a) interpret the integrant $ g(v) dg$ as a measure on the set
$G(\calA)f^{\lam}G(\calA)/G(\calA)$

\noindent
b) to show that the corresponding integral is convergent.

\medskip
\noindent
 Let us note that the expression $g(v)$ doesn't make
sense on $\GG$ but only on $\hGG$; more presicely, $g(v)$ makes
sense as $\VV$-valued function on the $\ZZ$-torsor over $\GG$
obtained from the $\bfK^*$-torsor $\hGG\to \GG$ by applying the
homomorphism $\bfK^*\to \ZZ$. Also the shift by $1\in \ZZ$
multiplies it by $q^{-h^{\vee}}$.

Denote by $\calT$ the $-h^{\vee}$-th power of this $\ZZ$-torsor.
Thus $g(v)$ makes sense as a $\VV$-valued function on $\calT$ so
that $1\in \ZZ$ acts by multiplication by $q$.
For a) we remind some basics of the theory of integration on varieites
over local fields such as $\bfE$.
Let  $Y$ be a smooth algebraic variety of $\bfE$ and let
$\Ome_Y$ denote the line bundle of highest forms on $Y$.
This line bundle defines a $\bfE^*$-torsor over $Y(\bfE)$.
By composing it with the valuation homomorphism $\bfE^*\to\ZZ$ we get a $\ZZ$-torsor
over $Y(\bfE)$. We denote by $\calT_Y$ the dual $\ZZ$-torsor.
Let $s$ be a locally constant function on $\calT_Y$ (with values in any vector space $V$) on which the shift by
$1\in\ZZ$ acts by means of multiplication by $q$.
Then it is well known that $s$ defines a $V$-valued measure on $Y(\bfE)$ for which
all open compact subsets are measurable. In particular, $\int_K s\in V$ makes sense where
$K$ is any open compact subset of $X(\bfE)$.
In our case it is easy to see that  we have an equality
$$
G(\calA)f^{\lam}G(\calA)/G(\calA)=G(\bfE[[f]])f^{\lam}G(\bfE[[f]]).
$$
The latter set is known to be the set of $\bfE$-points of a smooth variety $\Gr_G^{\lam}$
(the union of all $\Gr_G^{\lam}$'s is the {\it affine Grassmannian} $\Gr_G$ of $G$ - cf. \refs{grass}).
Moreover, the main property of the critical level will be the following:
\prop{tors-crit}
The restriction of $\calT$ to $G(\calA)f^{\lam}G(\calA)/G(\calA)$ is isomorphic
to $\calT_{\Gr_G^{\lam}}$.
\eprop
So we see that the part a) is achieved and we can consider
$ g(v) dg$ as a measure on the set
$G(\calA)f^{\lam}G(\calA)/G(\calA)$. If the variety $\Gr_G^{\lam}$ is proper (and thus the set $\Gr_G^{\lam}(\bfE)$ is compact) then the convergence of
the above integral would be obvious.
However, this happens only in some very limited number of cases, so theoretically some convergence issues
may arise. However,we show in \refs{grass} that the above integral is always
absolutely convergent.

\bigskip
\noindent
{\bf Remarks}.
\noindent

1. Note that the algebra $\CC[\Lam_+]$ is (abstractly) isomorphic
to the usual "finite" spherical Hecke algebra $\hsp_f$
 for the group $G(\bfK)$\footnote{Here the subscript "f" stands for "finite"}. In \refss{hom-finite} we make
this isomorphism canonical by looking at the action of $\hsp_{crit}$ in the
space of $G(\calA)$-invariants in some particular representation. However, apart from the case of
$G=GL(n)$ we can't compute this isomorphism explicitly.

\smallskip
\noindent
2. We expect that the algebra $\calR[\Lam_+]$ is dense in $\hsp_{crit}$ in some sense.
A precise conjecture is formulated in \refss{density}.
\ssec{variant}{A variant}Here we want to discuss what happens when the group $G(\calA)$ is replaced by
$G(\calO_{\bfF})$ (as was remarked above this case is somewhat easier but seems to be less relevant for global
applications discussed below).
In this case we should consider the functor
$$
\Coinv_c:\Rep_c(\hGG)\to \VV ect,
$$
sending every representation to its coinvariants with respect to $G(\calO_{\bfF})$.
We set
$$
\hspp_c=\End(\Coinv_c).
$$
As before we believe that the following is true:
\conj{go}
\begin{enumerate}
\item
The algebra $\hspp_c$ is commutative for every $c\in\CC$.
\item
One has $\hspp_c=\CC$ unless $q^{c+h^{\vee}}$ is a root of unity.
\end{enumerate}
\econj
Let us write $\hspp_{crit}$ instead of $\hspp_{-h^{\vee}}$. Then
We have the following result
\th{main'}There exists an embedding
$\CC[\Lam_+]\hookrightarrow \hspp_{crit}$.
\eth
The proof of \reft{main'} is basically a word-by-word repetition
of the proof of \reft{main} which is discussed in \refs{hecke} except that the last part
of the proof (\refss{end}) is not needed in the current case.

In fact, we believe that \reft{main'} may be strengthened in the following way.
\conj{density'}
The embedding in \reft{main'} is an isomorphism.
\econj
The reader should think of \refco{density'} as an analog of the density \refco{density}.

One can also formulate the following Conjecture which implies both \refco{go} and \refco{density'}.
Let $\overset{..}\bfH_c$ as before denote the modified double affine Hecke algebra as defined
in \cite{K2} and let $Z(\overset{..}\bfH_c)$ denote its center.
It follows easily from the results of \cite{K2} and \cite{GK1} that one has a natural map
$j_c:Z(\overset{..}\bfH_c)\to \hspp_c$.
\conj{}
The map $j_c$ is surjective for any $c$.
\econj
\ssec{global}{Global dreams}Here we would like to give some
(very speculative) motivation for considering the above Hecke algebras.
Let $S$ be a smooth surface over $\FF_q$ and
let $X\subset S$ be a smooth geometrically irreducible projective
curve in $S$. We assume for simplicity that
we also have a retraction $S_X\to X$ which is equal to identity on
$X$. Let $S_X$ be the formal neighbourhood of $X$
in $S$ and by $S_X^0$ the formal punctured neighborhood. We define in
\refs{glob}
 the set $\Bun_G(S_X^0)$ of isomorphism classes of
$G$-bundles on $S_X^0$ (over $\FF_q$) and a canonical $\Pic(X)(\FF_q)$-torsor ${\widehat\Bun}_G(S_X^0)$
over the set $\Bun_G(S_X^0)$. Using the degree map
$\Pic(X)\to \ZZ$ we obtain a $\ZZ$-torsor $\tbg(S_X^0)$ over
$\Bun_G(S_X^0)$.
Let $W$ be a finite collection of closed points of $X$ as a scheme over $\FF_q$.
Denote by $S_W$ the formal neighborhood of $W$ in $S$ and set
$S_W^0=S_W\cap S_X^0$. For every $w\in W$ we denote by $\bfK_w$ the local field
of $w$ in $X$ and by $\bfF_w$ the corresponding 2-dimensional
local field (the completion of $S_X^0$ along $\Spec \bfK_w$).
We let
$$
\bfF_W=\prod\limits_{w\in W}\bfF_w.
$$
Similarly let $\calA_w\subset \FF_w$ be the subring considered in \refss{kolcoa}  and let
$\calO_{S,w}$ denote the formal ring of $S$ at $w$. Clearly we have
$\calO_{S,w}=\calA_{w}\cap \calO_{\bfF_w}$ and we denote by $\calR_w$
the corresponding free commutative algebra whose generators are good curves in $\Spec (\calO_{S,w})$.
Let also
$$
\calA_W=\prod\limits_{w\in W}\calA_w.
$$
Let us also consider the group $\GG_W=G(\bfF_W)$. We also denote by $\tgg_W$ the corresponding
central extension of $\GG_W$ by means of $\ZZ$.
Let $\Bun_G(S_X^0,W)(\FF_q)$ denote the set of isomorphism
classes of the following data:

$\bullet$ A principal $G$-bundle $\calM$ on $S_X^0$ defined over
$\FF_q$

$\bullet$ A trivialization of $\calM$ on $S_W^0$.

We denote by $\tbg(S_X^0,W)$ the pull-back of the $\ZZ$-torsor
$\tbg(S_X^0)\to \Bun_G(S_X^0)$ under the natural map $\Bun_G(S_X^0,W)\to \Bun_G(S_X^0)$.
Note that the ring $\calA_W$ is exactly the ring of regular functions on $S_W^0$. Thus it
is clear that the group $G(\calA_W)$ acts on $\Bun_G(S_X^0,W)$ by changing the trivialization
of $\calM$ on $S_W^0$.

The following lemma is proved in \refs{glob}.
\lem{action-glob}
The action of $G(\calA_W)$ on $\Bun_G(S_X^0,W)$ extends canonically to an action
of $\tgg_W$ on $\tbg(S_X^0,W)$ in such a way that the central $\ZZ\subset \tgg_W$
acts on $\tbg(S_X^0,W)$ in the natural way.
\elem
We expect to be able to define
some ``correct'' pro-space of functions $\calF_W$ on $\tbg(S_X^0,W)$ which should be a representation
of the group $\tgg_W$. In particular, we should also have the pro-space  $\calF$
of functions on $\Bun_G(S_X^0)$ by taking $W$ to be the empty set. It can also be defined by
\eq{glob-inv}
\calF=\calF_W^{G(\calA_W)}.
\end{equation}
 For any $c\in\CC$ we should
also consider the spaces $\calF_{W,c}$ and $\calF_c$ of functions on which $1\in \ZZ$ acts
by means of multiplication by $q^c$. Generalizing \refe{glob-inv}, for every decomposition
$W=W'\coprod W''$ we should
have
\eq{glob-inv'}
\calF_{W',c}=\calF_{W,c}^{\GG(\calA_{W''})}.
\end{equation}
Now combining \refe{glob-inv'} with \reft{main} we see that for every $u\not\in W$ we get an action of
$\calR_u[\Lam_+]$ on $\calF_{W,-h^{\vee}}$. Moreover, for different choices of $u$ the corresponding
endomorphisms of $\calF_{W,-h^{\vee}}$ commute with each other. In other words, each
$\calF_{W,-h^{\vee}}$ is endowed with an action of the commutative algebra
$$
\bigotimes\limits_{u\not\in W} \calR_u[\Lam_+].
$$
One should think about this action as the "affine" analog of Hecke operators acting on the space of automorphic
forms for the group $G$ associated with the global field of rational functions on $X$. It would therefore
be very interesting to construct some examples of common eigen-functions of the above operators.
We believe that some generalizations of Eisenstein series considered in \cite{K1} should provide
such examples.
\ssec{}{Acknowledgements}
We would like to thank A.~Beilinson,
R.~Bezrukavnikov, V.~Drinfeld, D.~Gaitsgory, D.~Kaledin and S.~Kumar for helpful discussions
on the subject. The first author is partially supported by the NSF grant DMS-0300271.
\sec{grass}{Some results on singularities of affine Schubert varieties}
In this section we collect some fact about the canonical class of certain Schubert varieties that we shall
need in the future. All the results of this section follow easily from \cite{Fal} and \cite{BD}. In what follows we fix
a ground field $k$ and set $\calK=k((x))$ and $\calO=k[[x]]$.
\ssec{grass}{The affine Grassmannian}Let $G$ be a split semi-simple group over $k$ and set $\Gr_G=G(\calK)/G(\calO_K)$.
It is known that $\Gr_G$ has a natural structure of a proper ind-scheme over $k$. More precisely,
it is known that the orbits of the group $G(\calO)$ on $\Gr_G$ are parameterised by the elements of
$\Lam_+$. For each $\lam\in\Lam_+$ we shall denote by $\Gr_G^{\lam}$ the corresponding
orbit and by $\oGr_G^{\lam}$ its closure $\Gr_G$. The following theorem is proved in \cite{Fal}
(cf. also \cite{kum} and \cite{mat} for the corresponding result in characteristic 0).
\th{fal}
\begin{enumerate}
\item
Each $\oGr_G^{\lam}$ is a normal and Cohen-Macaulay projective variety over $k$.
\item
Each $\oGr_G^{\lam}$ has a resolution of singularities
\footnote{Of course, this statement is not empty only if char$k>0$.}
and for every such resolution $\pi^{\lam}:{\tGr_G^{\lam}}\to \oGr_G^{\lam}$
one has
$$
R\pi^{\lam}_*(\calO_{{\tGr_G^{\lam}}})=\calO_{\oGr_G^{\lam}}.
$$
(in other words $\oGr_G^{\lam}$ has {\it rational singularities}).
\end{enumerate}
\eth
\noindent
{\bf Remark.} In fact the above theorem is proved in \cite{Fal} in a slightly different context. Namely,
let $\Fl=G(\calK)/I$ denote the affine flag variety as in \cite{Fal}, where $I$ is an Iwahori
subgroup of $G(\calO)$ consisting of those elements of $G(\calO)$ whose reduction mod $t$ lies
in a fixed Borel subgroup of $G$. The $I$-orbits of on $\Fl$ are parametrised by the elements
of the affine Weyl group $W_{aff}$; for every $w\in W_{aff}$ we let $\Fl^w$ denote the
corresponding orbit and by $\oFl^w$ its closure in $\Fl$. The above theorem is proved in
\cite{Fal} for the varieties $\oFl^w$ rather than $\oGr_G^{\lam}$. However, it is easy to deduce
the statement for $\oGr^{\lam}$'s from the statement for $\oFl^w$'s. Namely, we have
the natural proper smooth projection $p:\Fl\to \Gr_G$. It is well-known that for every
$\lam\in\Lam_+$ there exist two elements $w_{\lam},w^{\lam}\in W_{aff}$ such that:
1) The projection $p$ maps $\oFl^{w_{\lam}}$ to $\oGr^{\lam}$ and the resulting map is proper and
birational.
2) $\oFl^{w^{\lam}}=p^{-1}(\oGr_G^{\lam})$.
It follows from 1) that any resolution of $\oFl^{w_{\lam}}$ is automatically a resolution of $\oGr_G^{\lam}$;
this shows that the latter variety has a resolution of singularities. Now given the existence of a resolution,
all the properties
claimed in \reft{fal} are properties of the singularities of $\oGr_G^{\lam}$; in particular, it is easy to see that
they are equivalent under smooth base change. Hence 2) shows that it is enough to check them for
$\oFl^{w^{\lam}}$.
\ssec{}{Line bundles on $\Gr_G$}It is well-known
(cf. \cite{BD} and \cite{Fal})
 that every finite-dimensional representation $V$ of
$G$ gives rise to a (determinant)
line bundle $\calL_V$ on $\Gr_G$. In particular, we let $\calL_{\grg}$ (here $\grg$ is the
Lie algebra of $G$) denote the line bundle corresponding to the adjoint representation of $G$.
We let $\calL_{crit}$ be the "critical line bundle" on $\Gr_G$; by the definition this is the (unique) square
root of $\calL_{\grg}^{-1}$.

In fact, it is well-known (cf. \cite{Fal}) that $\Pic(\Gr_G)\simeq \ZZ$. Under this identification
the line bundle $\calL_{\grg}$ is the $2h^{\vee}$-th power of a generator of $\ZZ$; we shall denote this
generator by $\calD$.

Let $\hatG(\calK)$ be the central extension of $G(\calK)$ by means of $k^*$ constructed in the same way as the one discussed
in \refss{notation} (with $k$ playing the role of $\bfK$ and $x$ playing the role of $t$).
Then $\calL_{crit}$ is $\hatG(\calK)$-equivariant in such a way that the central $k^*$ acts on the fibers of $\calL_{crit}$
by means of the character $a\mapsto a^{-h^{\vee}}$.
The following result from \cite{BD} is crucial for us:
\th{BD}
For every $\lam\in\Lam_+$ there is a canonical isomorphism
$$
\calL_{crit}|_{\Gr_G^{\lam}}\simeq \Ome_{\Gr_G^{\lam}}.
$$
(Here, as in the Introduction, $\Ome_{\Gr_G^{\lam}}$ denotes the bundle of highest forms
on $\Gr_G^{\lam}$.
\eth
The next result is an easy corollary of \reft{fal} and \reft{BD}.
\th{canonical}
\begin{enumerate}
\item
For every $\lam\in\Lam_+$ the variety $\oGr_G^{\lam}$ is Gorenstein. Moreover,
the canonical sheaf of $\oGr_G^{\lam}$ is isomorphic to $\calL_{crit}|_{\oGr_G^{\lam}}$.
Abusing the notation we shall denote this sheaf by $\Ome_{\oGr_G^{\lam}}$.
\item
For any $\lam\in\Lam_+$  let $\pi^{\lam}:{\tGr_G^{\lam}}\to \oGr_G^{\lam}$ be any resolution of
singularities. Then the identification between $(\pi^{\lam})^*\Ome_{\oGr_G^{\lam}}$ and
$\Ome_{{\tGr_G^{\lam}}}$ that one has at the generic point of ${\tGr_G^{\lam}}$ comes from an embedding
$$
(\pi^{\lam})^*\Ome_{\oGr_G^{\lam}}\hookrightarrow \Ome_{{\tGr_G^{\lam}}}.
$$
(In the case $\operatorname{char}k=0$ this implies that $\oGr_G^{\lam}$ has {\it canonical singularities}).
\end{enumerate}
\eth
\prf
Let us first prove (1). Let $\Ome_{\oGr_G^{\lam}}$ denote the Grothendieck-Serre dualizing complex
of $\oGr_G^{\lam}$ shifted by $-\dim \Gr_G^{\lam}$ (so that it coincides with $\Ome_{\Gr_G^{\lam}}$
on $\Gr_G^{\lam})$. According to \reft{fal} the variety $\oGr_G^{\lam}$ is Cohen-Macaulay.
Hence $\Ome_{\oGr_G^{\lam}}$ is in fact a sheaf (and not a complex of sheaves) which is automatically
reflexive. Equivalently, one can reformulate this as follows. Let $j^{\lam}:\Gr_G^{\lam}\hookrightarrow \oGr_G^{\lam}$
denote the natural open embedding. Since $\oGr_G^{\lam}$ is normal and
the complement to $\Gr_G^{\lam}$ in $\oGr_G^{\lam}$ has codimension
2 we have $\Ome_{\oGr_G^{\lam}}=j^{\lam}_*\Ome_{\Gr_G^{\lam}}$.
Now the normality of $\oGr_G^{\lam}$ also implies that
$$
\calL_{crit}|_{\oGr_G^{\lam}}=j^{\lam}_*(\calL_{crit}|_{\Gr_G^{\lam}})=j^{\lam}_*\Ome_{\Gr_G^{\lam}}=\Ome_{\oGr_G^{\lam}}.
$$
Let us now prove (2). According to \reft{fal} we have
$R\pi^{\lam}_*\calO_{\tGr_G^{\lam}}=\calO_{\oGr_G^{\lam}}$.
Applying the Grothendieck-Serre duality to both sides of this equality
we get $R\pi^{\lam}_*\Ome_{\tGr_G^{\lam}}=\Ome_{\oGr_G}^{\lam}$. By using the fact that
the functor $\pi^{\lam}_*$ is right adjoint to $(\pi^{\lam})^*$ we get the
map $(\pi^{\lam})^*\Ome_{\oGr_G}^{\lam}\to \Ome_{\tGr_G^{\lam}}$
which extends the natural
identification of these line bundles at the generic point of $\tGr_G^{\lam}$.
\epr
\ssec{integration}{Application to integration}
Let now $\bfE$ be a local non-archimedian field and let $Y$ be a smooth variety over $\bfE$.
Let $\calL$ be a line bundle on $Y$. Denote by $\calL_{\bfE^*}$ the corresponding $K^*$-torsor
over $Y(\bfE)$ and by $\calT_{\calL}$ the $\ZZ$-torsor obtained from $\calL_{\bfE^*}$ by pushing-forward
under the valuation map $\bfE^*\to\ZZ$. Set $\calS(Y,\calL)$ to be the space of
locally constant $\CC$-valued functions $s$ on $\calT_{\calL^{-1}}$ satisfying the following properties:

1) The element $1\in\ZZ$ acts on $s$ by means of
multiplication by $q$;

2) The support of $s$ is equal to the preimage of a compact subset of $Y(\bfE)$.

Note that every embedding $\calL_1\to\calL_2$ between two invertible sheaves on $Y$
gives rise to an embedding $\calS(Y,\calL_1)\to \calS(Y,\calL_2)$.
Similarly, for a pro-vector space $\VV$ over $\CC$ we may consider
the pro-vector space $\calS(Y,\calL)\ten \VV$.
Assume now that $Y$ is smooth and let $\calL=\Ome_Y$ (in this case we shall write $\calT_Y$ instead of
$\calT_{\Ome_Y^{-1}}$). Then it is well-known that we have a canonical
integration functional
$$
\int\limits_Y :\calS(Y,\Ome_Y)\to \CC.
$$
Thus for a pro-vector space $\VV$ we also get the integration map
$$
\int\limits_Y \calS(Y,\Ome_Y)\ten \VV\to\VV.
$$
Assume now that $Y$ is not necessarily smooth but satisfies the following conditions:
1) $Y$ is Gorenstein (we denote by $\Ome_Y$ the corresponding line bundle);
2) $Y$ has canonical singularities, i.e. there exists a resolution of singularities
$\pi:\tilY\to Y$ and an embedding $\pi^*\Ome_Y\to\Ome_{\tilY}$ extending the
natural identification of these line bundles at the generic point of $\tilY$.
Then we claim that the integration map $\int_Y$ as above is still well-defined.
Indeed, it follows from condition 2) that every $s\in\calS(Y,\Ome_Y)$ gives rise
to some $\tils\in\calS(\tilY,\Ome_{\tilY})$ and we define
$$
\int\limits_Y s=\int\limits_{\tilY}\tils.
$$
It is clear that $\int\limits_{\tilY}\tils$ does not depend on a choice of a
resolution $\tilY$ since in fact this integral is equal to the integral of $s$ over the smooth part of
$Y$ (the existence of $\tilY$ with the above properties shows that this integral is absolutely
convergent).
So, for any pro-vector space and any $\lam\in\Lam_+$ we have a well-defined
integration map
\eq{integration}
\int\limits_{\oGr_G^{\lam}}:\calS(\oGr^{\lam},\calL_{crit})\ten\VV\to \VV.
\end{equation}
Sometimes we shall write $\int\limits_{\Gr_G^{\lam}}$ instead of $\int\limits_{\oGr_G^{\lam}}$.
\sec{hecke}{The Hecke algebras}
In this Section we  prove \reft{main} and discuss some corollaries of it.
First, let us discuss the functor of invariants in more detail.
\ssec{invc}{The functor $\Inv_c$}Since we are dealing with pro-vector spaces let us first explain what we mean by
the functor $\Inv_c$ of $G(\calA)$-invariants on the category $\Rep_c(\hGG)$. In \cite{GK1} the authors explain
the for any $\VV\in\Rep(\hGG)$ and any subgroup $H\in\hGG$ one can define a functor from $\VV ect$ to $Vect$
which is supposed to be represented by $\VV^H$. However, this functor is representable only if certain conditions
are satisfied. In this sub-section we want to check that these conditions are satisfied in the case of the subgroup
$G(\calA)\subset \hGG$. In fact this is not necessary for the rest of this section. Namely in the next
sub-section
we are going to construct some endomorphisms of $\VV^{G(\calA)}$ (for every $\VV$). For this purpose we shall only
need the functorial definition $V^{G(\calA)}$. However, the fact that $V^{G(\calA)}$ actually exists as a pro-vector
space shows the existence of the functor $\Inv_c:\Rep_c(\hGG)\to \VV ect$
and  it still nice to know that the endomorphisms that we are going to construct
will actually be endomorphisms of this functor.

Let us recall some definitions from \cite{GK1}.

Let us denote by $Set_0$ the category of finite sets. Let also
$$
\bfS et=\Ind(\Pro(Set_0));\qquad \SS et=\Ind(\Pro(\bfS et)).
$$

We have the natural functor $T:Set_0\to \bfS et$.
We shall need the functor
$$
\TT:\bfS et=\Ind(\Pro(Set_0))\to \SS et=\Ind(\Pro(\Ind(\Pro(Set_0)))
$$
which is defined by first applying $T$ to the "inner" $Set_0$ and then applying Ind(Pro) to both sides.

Let $\HH$ be an object of $\SS et$. In \cite{GK1} the notion of action of $\HH$ on
a pro-vector space $\VV$ is defined. Also if $\HH$ is a group-like object
of $\SS et$ then the notion of a representation of $\HH$ on a pro-vector space is defined
and we also have a well defined category $\Rep(\HH)$ of pro-representations of $\HH$.
One has the functor $\triv:\VV ect\to \Rep(\HH)$ corresponding to "trivial" representations.
It is shown in \cite{GK1} (cf. Proposition 2.10) that this functor always admits a left adjoint
(which should be called the functor of coinvariants of $\HH$). In addition, triv also admits
a right adjoint functor if the following technical condition is satisfied:

\medskip
\noindent
(**) As an object of $\SS et$ the group $\HH$ can be represented as a limit
$"\underset{\to}\lim" X_k$\footnote{Here $"\underset{\to}\lim"$ means that the limit is taken as an object
of the corresponding Ind-category} where each $X_k\in \Pro(\bfS et)$ is weakly strict (we refer the reader to
\cite{GK1} for the definition of this notion).

\medskip
\noindent
In our case let $\HH=G(\calA)$.
The ring $\calA$ naturally gives rise to an object of $\bfS et$ since
$$
\calA=\underset{\to}\lim \ t^{-k}\calO_{\bfK}[[t]]
$$
and the ring $\calO_{\bfK}[[t]]$ is clearly a projective limit of finite sets.
However, we are going to regard the ring $\calA$ as an object in $\SS et$ by applying the functor
$\TT$ to the above construction.
For any affine algebraic variety $X$ over $\calA$ we may also regard $X(\calA)$ as an object of
$\SS et$ by embedding it into $\calA^n$. Moreover, $X(\calA)$ always satisfies the condition (**)
(it is enough to check this for $\calA$ itself and in this case this is obvious).

Thus the group $G(\calA)$ may be regarded as a group-like object of the category $\SS et$ which
satisfies (**). Hence the functor of $G(\calA)$-invariants is well-defined on the category $\Rep(G(\calA))$
according to {\it loc. cit.} Also, the embedding $G(\calA)\to \hGG$ is an embedding of group-like
objects of $\SS et$. Hence the restriction functor $\Rep_c(\hGG)\to \Rep(G(\calA))$ is well-defined.
By composing these two functors we get the functor
$$
\Inv_c:\Rep_c(\hGG)\to \VV
$$
of $G(\calA)$-invariants on the category $\Rep_c(\hGG)$.

\ssec{}{The construction of $\alp_C$}We want to define the homomorphism $\alp_C$ by means of the
formula \refe{alpC}. First of all, without loss of generality we may
assume that $f=x$.
Consider the ring $\calB=\FF_q[[x,t]]$ and let $\bfF_{\calB}$ denote its field of fractions.
It is easy to see that the natural embedding
$\calA=\FF_q[[x]]((t))\hookrightarrow \FF_q((t))[[x]]=\bfE[[x]]$ induces an isomorphism
\eq{agrl}
G(\calA)x^{\lam}G(\calA)/G(\calA){\widetilde \to} G(\bfE[[x]])x^{\lam}G(\bfE[[x]])=\Gr_G^{\lam}(\bfE).
\end{equation}
We need now the following result which generalizes slightly \refp{tors-crit}.
It follows from \refe{agrl} that we have the natural identification
\eq{agrl'}
 \bigcup\limits_{\mu\leq \lam}G(\calA)x^{\mu}G(\calA)/G(\calA)=\oGr_G^{\lam}(\bfE).
\end{equation}
The left hand side has a natural $\bfK^*$-torsor over it coming from the central extension $\hGG$.
As before we denote by $\calT$ the $-h^{\vee}$-th power of the corresponding $\ZZ$-torsor.
\prop{tors-crit'}
\begin{enumerate}
\item
Under the identification \refe{agrl'}
the $\ZZ$-torsor $\calT$ defined above goes over to $\calT_{\oGr_G^{\lam}}$.
\item
Let now $\VV$ be a $\hGG$-representation with the critical central charge.
Then the formula $v\mapsto g(v)$ for
$$
g\in \bigcup\limits_{\mu\leq \lam}G(\calA)x^{\mu}G(\calA)/G(\calA)
$$
defines a linear map $\VV^{G(\calA)}\to \calS(\Gr_G^{\lam},\Ome_{\Gr_G^{\lam}})\ten \VV^{G(\calA)}$.
\end{enumerate}
\eprop
First of all, let us explain why \refp{tors-crit'} implies the construction of $\alp_C$.
 By composing the map $\VV^{G(\calA)}\to \calS(\Gr_G^{\lam},\Ome_{\oGr_G^{\lam}})\ten \VV^{G(\calA)}$
 defined in (2) of \refp{tors-crit'}
it with the integration map $\int_{\Gr_G^{\lam}}$ as in  \refss{integration} we get a well defined
linear map $\VV^{G(\calA)}\to \VV^{G(\calA)}$ which is clearly functorial in $\VV$. In other
words we get an element of $\hsp_{crit}$ which by definition is equal to $\alp_C(\lam)$.

For the proof of  \refp{tors-crit'} we need the following
auxiliary result.
By the definition the field $\bfF_{\calB}$ is embedded naturally into the field $\bfF=\FF_q((x))((t))$ as well as into the
field $\FF_q((t))((x))$. Hence the group $G(\bfF_{\calB})$ acquires two natural central extensions:
one by means of $\bfK^*=\FF_q((x))^*$ and the other by means of
$\bfE^*=\FF_q((t))^*$.  By applying the valuation homomorphisms
$\bfK^*\to \ZZ$ and $\bfE^*\to \ZZ$ we get two central extensions of $G(\bfF_{\calB})$ by means
of $\ZZ$. Clearly, these two central extensions are interchanged by the automorphism
of $G(\bfF_{\calB})$ induced by the automorphism of $\bfF_{\calB}$ which interchanges $x$ and $t$.
\lem{exten}
The above two central extensions are opposite to each other.
\elem
\prf
Recall that for a field $R$ we have a natural central extension
\eq{kdva}
1\to K_2(R)\to \oG(R)\to G(R)\to 1
\end{equation}
where $K_2(R)$ denotes the corresponding Milnor $K$-group.
In the case when $R=k((t))$ (where $k$ is any field) we have the tame symbol homomorphism
$K_2(R)\to k^*$ defined as follows. Recall that $K_2(R)$ is a quotient of $\Lam^2(R^*)$ by means
of all the expressions of the form $f\wedge (1-f)$.
Thus to define the above homomorphism it is enough to define the tame symbol $(f,g)\in k^*$ for
any $f,g\in R^*$. Let $v: R^*\to \ZZ$ be the valuation homomorphism. One defines
$$
(f,g)=(-1)^{v(f)v(g)}\frac{f^{v(g)}}{g^{v(f)}}(0).
$$
The fact that $(\cdot,\cdot)$ descends to $K_2(R)$ is well-known (and straightforward).

It is also well-known (cf. \cite{BD}) that by pushing forward the extension \refe{kdva}
 with respect to the tame symbol
we get the central extension
$$
1\to k^*\to\tilG(R)\to G(R)\to 1
$$
that was discussed above.

The same construction also applies if $R$ is a subfield of $k((t))$.

Let now $R$ be the field of fractions of the ring $k[[x,t]]$. Then by regarding $R$ as a subfield of both
$k((x))((t))$ and $k((t))((x))$ we get two tame symbols
$$
(\cdot,\cdot)_t:K_2(R)\to k((x))^*\quad\text{and}\quad (\cdot,\cdot)_x:K_2(R)\to k((t))^*
$$
and two valuation homomorphisms
$$
v_t:k((t))^*\to \ZZ\quad\text{and}\quad v_x:k((x))^*\to \ZZ.
$$
Thus \refl{exten} follows from the following equality whose verification we leave to the reader:
for any $f,g\in R^*$ we have
$$
v_x((f,g)_t)=-v_t((f,g)_x).
$$
\epr
It is now clear that the assertion (1) of \refp{tors-crit'} follows from \refl{exten} combined with \reft{BD}.
From this the assertion (2) is straightforward.

\medskip
\noindent
It is easy to see
that the homomorphisms $\alp_C$ defined above together
define a homomorphism $\calR[\Lam_+]\to \hsp_{crit}$.
Indeed, for this we have to verify the following two statements:

\medskip
\noindent
1) For any $C\in\calC$ and $\lam,\mu\in\Lam_+$ we have
\eq{composition}
\alp_C(\lam)\cdot\alp_C(\mu)=\alp_C(\lam+\mu).
\end{equation}

\noindent
2) For any $C,C'\in\calC$ and $\lam,\mu\in\Lam_+$ the elements $\alp_C(\lam)$ and $\alp_{C'}(\mu)$
commute with each other.

First of all, if $C\neq C'$ then statement 2) is obvious: in this case both $\alp_C(\lam)\cdot\alp_{C'}(\mu)$ and
$\alp_{C'}(\mu)\cdot\alp_C(\lam)$ are given by integrals of the same measure over $\Gr^{\lam}(\bfE)\x \Gr^{\mu}(\bfE)$.
Hence it is enough to check 1) (note that \refe{composition} implies that $\alp_C(\lam)$ and $\alp_C(\mu)$ commute).
For any field $k$ consider the ind-scheme $\Gr_G *\Gr_G$ whose set of $k$ points is equal
$G(\calK)\underset{G(\calO)}\x G(\calK)/G(\calO)$ here as before we denote $\calK=k((x))$ and $\calO=k[[x]]$).
We have the natural maps
$$
p_1:G(\calK)\x\Gr_G\to \Gr\quad p_2:G(\calK)\x\Gr_G\to \Gr_G\quad\text{and}\quad
m:\Gr_G*\Gr_G\to \Gr_G
$$
which are defined by
$$
p_1(g_1,g_2\mod G(\calO))=g_1;\quad p_2(g_1,g_2\mod G(\calO))=g_2\mod G(\calO)
$$
and
$$
m(g_1,g_2 \mod G(\calO))=g_1g_2 \mod G(\calO).
$$
We also have the natural projection $\pi:G(\calK)\to \Gr_G$.
For any $\lam,\mu\in\Lam_+$ let us set $\Gr_G^\lam *\Gr_G^\mu$
to be the image of $(\pi\circ p_1)^{-1}(\Gr_G^\lam)\x p_2^{-1}(\Gr_G^\mu)\subset G(\calK)\x \Gr_G$ in
$\Gr_G*\Gr_G$. We also denote its closure by $\oGr_G^\lam*\oGr_G^\mu$.
The following result is well-known:
\lem{}
The map $m$ maps $\Gr_G^\lam*\Gr_G^\mu$ to $\oGr_G^{\lam+\mu}$.
\elem
Hence we also have the map $m:\oGr_G^\lam*\oGr_G^\mu\to\oGr_G^{\lam+\mu}$ which is proper and birational.
Thus (using the same notations as in \refe{alpC} we get
$$
\alp_C^{\VV}(\lam)\cdot\alp_C^{\VV}(\mu)(v)=\int\limits_{g_1\in G(\calA)f^{\lam}G(\calA)/G(\calA)}
g_1(\int\limits_{g_2\in G(\calA)f^{\mu}G(\calA)/G(\calA)} g_2(v) dg_2)dg_1=
$$
$$
\int\limits_{g\in \Gr_G^\lam(\bfE)*\Gr_G^\mu(\bfE)} g(v) dg=
\int\limits_{g\in \Gr_G^{\lam+\mu}(\bfE)} g(v) dg=\alp_C^{\VV}(\lam+\mu)(v).
$$
\medskip
\noindent
It remains to show that this homomorphism is
an embedding. Before we proceed with this let us look at one example of the action of $\calR[\Lam_+]$
on $\VV^{G(\calA)}$ for some particular representation $\VV$.
\ssec{hom-finite}{The homomorphism $\hsp_{crit}\to\hsp_f$}
Let $\hsp_f$ denote the usual spherical Hecke algebra associated with the group $G$ and the local field
$\bfK$. This algebra can be defined as follows. Let $\calS(G(\calO_{\bfK})\backslash G(\bfK))$ be the space
of $\CC$-valued functions on $G(\calO_{\bfK})\backslash G(\bfK)$ with finite support. This is a representation of the
group $G(\bfK)$ (acting on the right).
Then
\eq{hecke-finite}
\hsp_f=\End_{G(\bfK)}\calS(G(\calO_{\bfK})\backslash G(\bfK)).
\end{equation}
Equivalently, $\hsp_f$ is the algebra of endomorphisms of the functor of $G(\calO_{\bfK})$-invariants
on the category $\Rep G(\bfK)$ of smooth representations of $G(\bfK)$.

Let $G^{\vee}$ denote the Langlands dual group of $G$ (considered as an algebraic group over $\CC$).
Thus $\Lam$ is the weight lattice of $G^{\vee}$ and $\Lam_+$ is the set of dominant weight of $G^{\vee}$.
Recall that $\hsp_f$ is canonically isomorphic to the complexified Grothendieck group
of the category of finite-dimensional
representations of $G^{\vee}$. In particular, $\hsp_f$ has a basis
$\{ A_{\lam}\}_{\lam\in\Lam_+}$ where each $A_{\lam}$ corresponds to the irreducible representation of
$G^{\vee}$ with highest weight $\lam$.

As is well-known the algebra $\hsp_f$ is also isomorphic to the algebra of
$G(\calO_{\bfK})$-biinvariant functions on $G(\bfK)$ (the algebra structure is with
respect to convolution; here we choose a Haar measure on $G(\bfK)$ which
is characterized by the property that the volume of $G(\calO_{\bfK})$ is equal to 1).
Since we have the natural identification
$$
G(\calO_{\bfK})\backslash G(\bfK)/G(\calO_{\bfK})=\Lam_+,
$$
it follows that
$\hsp_f$ has another basis $T_\lam$ where each $T_\lam$ is equal to the characteristic
function of the corresponding double coset multiplied
by $(-q)^{\la \lam,\rho^{\vee}\ra}$
\footnote{Here $\rho^{\vee}$ is the half-sum of all the positive roots of $G$}.
It is also well known that if we let $\rho^{\vee}$ denote the half-sum of the positive roots of $G$ then we have
$$
A_\lam=(-q^{1/2})^{\la\lam,\rho^{\vee}\ra}T_\lam+\text{linear combination of $T_\mu$ with $\mu<\lam$}
\footnote{Here we say that $\mu<\lam$ if the difference $\lam-\mu$ is a sum of positive
coroots}.
$$
Let now $\Gam_1\subset G(\calO_{\bfF})$ denote the first congruence subgroup (i.e. $\Gam_1$ is the kernel
of the natural "reduction mod $t$" homomorphism $G(\calO_{\bfF})\to G(\bfK)$).
Let $\calS(\hGG/\Gam_1)$ denote the Schwartz space
\footnote{Recall that  $\calS(\hGG/\Gam_1)$ is actually a pro-vector space}
of $\hGG/\Gam_1$ in the sense of \cite{GK1} (equivalently,
$\calS(\hGG/\Gam_1)$ is the space of $\Gam_1$-coinvariants in the Schwartz space $\calS(\hGG)$ introduced in
\cite{GK1}). For $c\in \CC$ we also denote by $\calS_c(\hGG/\Gam_1)$ the coinvariants of $\bfK^*$ on $\calS(\hGG/\Gam_1)$
with respect to the character defined by $c$. This is a representation of $\GG\x G(\bfK)$ where the first factor
acts on the left and the second factor acts on the right.
In particular, we may consider the space $\calS_c(\hGG/\Gam_1)^{G(\calA)}$
of $G(\calA)$-invariants with respect to the left action. It has a natural action of $G(\bfK)$.
We have the natural embedding $\phi:G(\bfK)\to
\Gam_1 \backslash \GG$ which identifies the left hand side with
$\Gam_1 \backslash G(\calO_{\bfK})$. Since the central extension splits canonically over $G(\calO_{\bfK})$ we get a well-defined
restriction map
$$
\phi^*_c:\calS_c(\Gam_1  \backslash \hGG)\to \calS(G(\bfK)).
$$
Since the image of $G(\calA)\cap G(\calO_{\bfK})$ is equal to $G(\calO_{\bfK})$ it follows that
$\phi^*_c$ gives rise to the map
$$
\phi^*_{c,\calA}:\calS_c(\Gam_1  \backslash \hGG)^{G(\calA)}\to \calS(G(\calO_{\bfK})\backslash G(\bfK))
$$
which is a homomorphism of right $G(\bfK)$-modules.
\lem{rest-iso}
The map $\phi^*_{c,\calA}$ defined above is an isomorphism.
\elem
\prf
The proof follows immediately from the following result.
\lem{gago}
For any connected linear algebraic group $G$ one has $\GG=G(\calA)\cdot G(\calO_{\bfF})$.
\elem
\prf
Consider as before the scheme $\Del=\Spec(\FF_q[[x,t]])$ and let $0$ denote its unique
closed point. Let also $\Del^0$ denote the complement of $0$ and denote by $j:\Del^0\to \Del$ the natural
embedding. Let $U=\Spec \calO_{\bfF}$ and
$V=\Spec\calA$. Then both $U$ and $V$ are open subsets of $\Del^0$ and we have
$$
\Del^0=U\cup V;\quad U\cap V=\Spec\bfF.
$$
We claim now that any $g\in\GG$ defines a $G$-bundle $\calG_g$ on $\Del^0$ in the following way:
consider the trivial $G$-bundles on $U$ and $V$ respectively and use $g$ to glue them
on $U\cap V=\Spec\bfF$. Note that$\calG_g$ constructed in this way is automatically trivialized on
$U$ and $V$; changing these trivializations amounts to multiplying $g$ on the left by elements of
$G(\calO_{\bfF})$ and on the right by elements of $G(\calA)$.
It is clear that the isomorphism class of $\calG_g$
is determined by the class of $g$ in $G(\calO_{\bfF})\backslash \GG/G(\calA))$.
Thus to prove \refl{gago} it is enough to show that $\calG_g$ is isomorphic to the trivial
$G$-bundle. We claim that every $G$-bundle $\calG$ on $\Del^0$ is trivial.
 First of all observe that $\calG$ extends canonically (and uniquely) to $\Del$.
Indeed, if $\calM$ is a vector bundle on $\Del^0$ then it is known that $j_*\calM$ gives its unique
extension to a vector bundle on $\Del$. Now a $G$-bundle on $\Del^0$ is the same
as a tensor functor from the category of finite-dimensional $G$-modules to the category of
vector bundles on $\Del^0$. By applying $j_*$ to every such functor we get a tensor functor
from the category of finite-dimensional $G$-modules to the category of
vector bundles on $\Del$ which is the same as a $G$-bundle on $\Del$.
Now it is enough to show that every $G$-bundle on $\Del$ is trivial. This is obvious when $\FF_q$ is replaced
by ${\overline \FF_q}$. Hence over $\FF_q$ such $G$-bundles are classified
$H^1(\operatorname{Gal}({\overline\FF}_q/\FF_q),G_{ad})=0$
where $G_{ad}$ denotes the adjoint group of $G$. However, this group is known to be trivial
\footnote{Here we use in an essential way that $G$ is connected} which finishes the proof.
\epr
\epr

Let now $h\in \hsp_c$. By the definition it gives rise to an endomorphism of $\calS_c(\hGG/\Gam_1)^{G(\calA)}$
which commutes with every endomorphism of this space coming from an endomorphism of
$\calS_c(\hGG/\Gam_1)$ as a $\hGG$-representation. In particular, it commutes with the right translations
by elements of $G(\bfK)$. Hence it follows from \refl{rest-iso} that $h$ gives rise to an endomorphism of
$\calS(G(\calO_{\bfK})\backslash G(\bfK))$ which commutes with right translations by elements of $G(\bfK)$.
By \refe{hecke-finite} this is the same as an element of $\hsp_f$. In this way we get
a homomorphism
$$
\iota_c:\hsp_c\to \hsp_f.
$$
When $c=-h^{\vee}$  we shall denote $\iota_c$ by $\iota_{crit}$.
\lem{homo}
For every $C\in\calC$ we have
$$
\iota_{crit}\circ\alp_C(\lam)=(-q^{\frac{1}{2}})^{\la \lam,\rho^{\vee}\ra}
A_{\lam}+\text{linear combination of $A_{\mu}$'s with $\mu<\lam$}=
$$
$$
T_{\lam}+\text{linear combination of $T_{\mu}$'s with $\mu<\lam$}.
$$
In particular, $\iota_{crit}\circ\alp_C$ is an isomorphism.
\elem
\prf
First of all we want to show that $\iota_{crit}\circ\alp_C(\lam)$ lies in the span
of $A_\mu$ with $\mu\leq \lam$. This is equivalent to saying that
$\iota_{crit}\circ\alp_C(\lam)$ lies in the span of $T_\mu$ with $\mu\leq \lam$.
It is clear that it is enough to assume that $C$ is given by the equation
$x=0$.

To prove this assertion we must look at the support of $\iota_{crit}\circ\alp_C(\lam)$
considered as a $G(\calO_{\bfK})$-biinvariant function on $G(\bfK)$.
In other words, we arrive to the following "set-theoretic" problem.
We have the identification
$$
\Gam_1\backslash \GG/G(\calA)=G(\bfK)/G(\calO_{\bfK})=\Gr_G(\FF_q).
$$
Thus we obtain the map $\del:\GG/G(\calA)\to \Gr_G(\FF_q)$. We need to show that the image
of $G(\calA)x^{\lam}G(\calA)/G(\calA)$ under $\del$ is contained in
$\oGr_G^{\lam}(\FF_q)$. Recall that we denote by $\bfE$ the field $\FF_q((t))$ and
that we have the natural identification
$$
G(\calA)x^{\lam}G(\calA)/G(\calA)=\Gr_G^{\lam}(\bfE).
$$
Since $\Gr_G$ is ind-proper it follows that we have $\Gr_G(\bfE)=\Gr_G(\calO_{\bfE})$ (here
$\calO_{\bfE}=\FF_q[[t]]$ thus we have the well-defined "reduction mod $t$" map
$$
G(\calA)x^{\lam}G(\calA)/G(\calA)=\Gr_G^{\lam}(\bfE)\to \Gr_G(\FF_q).
$$
It is easy to see that this map actually coincides with the restriction of $\del$ to
$G(\calA)x^{\lam}G(\calA)/G(\calA)$. On the other hand, since $\oGr_G^{\lam}$ is proper,
the above map actually lands in $\oGr^{\lam}(\FF_q)$ which is what we had to prove.

To compute the coefficient of $T_{\lam}$ in $\iota_{crit}\circ\alp_C(\lam)$ (which is the same
as the coefficient of $A_\lam$ up to the factor of $(-q^{1/2})^{\lam,\rho^{\vee}}$) we need
to study the fibers of the map
$$
G(\calA)x^{\lam}G(\calA)/G(\calA)\cap \del^{-1}(\Gr_G^{\lam}(\FF_q))\to \Gr_G^{\lam}(\FF_q).
$$
Since $\Gr_G^{\lam}(\FF_q)$ is smooth, it follows from Hensel's lemma that each fiber is isomorphic
to $(t\calO_{\bfE})^{\dim \Gr_G^{\lam}}$ and its volume with respect to a measure coming from a differential
form on $\Gr_G^{\lam}$ defined over $\calO_{\bfE}$ is equal to 1. This finishes the proof.
\epr
\noindent
{\bf Remark.} Let $G=PGL(n)$. Then \refl{homo} describes the map $\iota_{crit}\circ\alp_C$
completely.
Indeed, since the group $PGL(n)$ is of adjoint type, the semigroup $\Lam_+$ is generated by
the fundamental coweights $\ome_1,...,\ome_{n-1}$ and thus it is enough to describe
all the $\iota_{crit}(\ome_i)\circ\alp_C$. Also in this case for each
$i=1,...,n-1$ the set
$$
\{ \mu\in\Lam_+|\ \mu<\ome_i\}
$$
is empty. Hence $\iota_{crit}\circ\alp_C(\ome_i)=A_{\ome_i}$.
We don't know how to describe $\iota_{crit}$ for other groups.
\ssec{end}{End of the proof}
We now want to show that the above homomorphism
$\calR[\Lam_+]\to \hsp_{crit}$ is an embedding.
To do that it is sufficient to check the following: let $C_1,...,C_k$ be pairwise distinct
elements of $\calC$ and let $\mu_1,...,\mu_k\in \Lam_+$. Consider the elements
\eq{elements}
\alp_{C_1}(\mu_1)\cdot ...\cdot\alp_{C_k}(\mu_k)\in\hsp.
\end{equation}
 We need to show that all these elements
 are linearly independent.

For this we want to see how these elements
act on the pro-space of $G(\calA)$-invariants for representations $\VV_n=\calS_{crit}(\hGG/\Gam_n)$ where $\Gam_n$ denotes
the $n$-th congruence subgroup (by definition $\Gam_n$ is the kernel of the natural
homomorphism $G(\calO_{\bfF})\to G(\calO_{\bfF}/\grm_{\bfF}^n)$).

Let $G_n=G(\calO_{\bfF}/\grm_{\bfF}^n)$ and let also $G_{n,\calO}=G(\calO_{\bfK}[[t]]/t^n)$.
Then $G_n$ is a locally compact totally disconnected group and $G_{n,\calO}$ is an open
compact subgroup of it. Also
we have the natural identification $\VV_n^{G(\calA)}=G_{n,\calO}\backslash G_n$. Thus
the image of $\hsp_{crit}$ in $\End(\VV_n^{G(\calA)})$ embeds into the Hecke
algebra $H(G_n,G_{n,\calO})$. In particular, we may think about the elements
$\alp_{C_1}^{\VV_n}(\mu_1),...,\alp_{C_k}^{\VV_n}(\mu_k)$ as $G_{n,\calO}$-biinvariant functions
on $G_n$.
Thus it is enough to show that for any finite collection of elements of the type \refe{elements} their
images in the images $H(G_n,G_{n,\calO})$ are linearly independent functions for $n$ large enough.

Let $f_1,...,f_k$ be the equations of the curves $C_1,...,C_k$ and
set
$$
X_{n,C_1,...,C_k}^{\lam_1,...,\lam_k}=G_{n,\calO}f_1^{\lam_1}...f_k^{\lam_k}G_{n,\calO}\subset G_n
\footnote{We shall always assume that $\lam_j\neq 0$ unless $k=0$}.
$$

\lem{cosets}
\begin{enumerate}
\item
The element $\alp_{C}^{\VV_n}(\lam)$ is supported on
$$
\cup_{\mu\leq \lam}X_{n,C}^\mu
$$
and doesn't vanish on $X_{n,C}^\lam$.
\item
Assume that the curves $C_1,...,C_k$ are distinct. Then for any collection
$((C_1',\mu_1),...,(C_l',\mu_l))$ different from $((C_1,\lam_1),...,(C_k,\lam_k))$ and for $n$ large enough
the set
$X_{n,C_1}^{\lam_1}\cdot ...\cdot X_{n,C_k}^{\lam_k}$
does not contain  $X_{n,C_1',...,C_l'}^{\mu_1,...,\mu_l}$.
\item
Fix a finite set of pairwise distinct collections $(C_1,\lam_1),...,(C_k,\lam_k))$. The for $n$ large enough
all the corresponding cosets $X_{n,C_1,...,C_k}^{\lam_1,...,\lam_k}$ are pairwise distinct.
\end{enumerate}
\elem
\prf
The first assertion is proved in exactly the same way as \refl{homo}. Let us prove (2) and (3).
Let us "take the limit" $n\to\infty$. In other words, let us recall the notation
$\calB=\calO_{\bfF}\cap \calA=\FF_q[[x,t]]$. For $(C_1,\lam_1),...,(C_k,\lam_k))$
as above set
$$
X_{C_1,...,C_k}^{\lam_1,...,\lam_k}=G_{\calB}f_1^{\lam_1}...f_k^{\lam_k}G_{\calB}\subset G(\calO_{\bfF}).
$$
Then we have
$$
X_{C_1,...,C_k}^{\lam_1,...,\lam_k}=\underset{\leftarrow}\lim \ X_{n,C_1,...,C_k}^{\lam_1,...,\lam_k}
$$
and hence
the assertions (2) and (3) of \refl{cosets} follow from the following two statements:

\medskip
a) Assume that the curves $C_1,...,C_k$ are distinct. Then for any
$$
((C_1',\mu_1),...,(C_l',\mu_l))\neq ((C_1,\lam_1),...,(C_k,\lam_k)),
$$
the set $X_{C_1}^{\lam_1}\cdot ...\cdot X_{C_k}^{\lam_k}$
does not contain  $X_{C_1',...,C_l'}^{\mu_1,...,\mu_l}$.

\smallskip
b) The cosets $X_{C_1,...,C_k}^{\lam_1,...,\lam_k}$ are  distinct for distinct collections
$((C_1,\lam_1),...,(C_k,\lam_k))$.

\medskip
\noindent
To prove a) and b) let us give a geometric interpretation of these statements.
Recall that we denote $\Del=\Spec \calB$. Then it is easy to see that the set
$G(\calB)\backslash G(\calO_{\bfF})/G(\calB)$ naturally parametrises the following
data:

\medskip
1) Two principal $G$-bundles $\calF_1,\calF_2$ on $\Del$;

\smallskip
2) An identification $\kap$ between $\calF_1$ and $\calF_2$ defined outside the union
of all the curves $C$ in $\Del$ which are different from the curve $t=0$ .

\medskip
\noindent
Under this identification, the coset $X_{C_1,...,C_k}^{\lam_1,...,\lam_k}$ lies in the subset
$Y_{C_1,...,C_k}^{\lam_1,...,\lam_k}$ of
all triples $(\calF_1,\calF_2,\kap)$ as above  such that
$\kap$ is well-defined outside the curves $C_1,...,C_k$ and at the generic point of every $C_j$
the singularity of $\kap$ is of type $\lam_j$ (note, however, that in general
$X_{C_1,...,C_k}^{\lam_1,...,\lam_k}\neq Y_{C_1,...,C_k}^{\lam_1,...,\lam_k}$
unless $k=1$). This immediately proves b) since the sets $Y_{C_1,...,C_k}^{\lam_1,...,\lam_k}$ are clearly
disjoint.
To prove a) it is enough to observe that  we have
$$
(\calF_1,\calF_2,\kap)\in G(\calB)\backslash X_{C_1}^{\lam_1}\cdot ...\cdot X_{C_k}^{\lam_k}/G(\calB)
$$
if and only if there exists a chain
$$
(\calG_0,\calG_1,\kap_1),\ (\calG_1,\calG_2,\kap_2),...,
(\calG_{k-1},\calG_k,\kap_k)\in G(\calB)\backslash G(\calO_{\bfF})/G(\calB)
$$
such that

(i) $\calG_0=\calF_1$, $\calG_k=\calF_2$;

(ii)$\kap$ is equal to the composition of all the $\kap_j$;

(iii) $(\calG_{j-1},\calG_j,\kap_j)\in G(\calB)\backslash X_{C_j}^{\lam_j}/G(\calB)$ for every $j=1,...,k$.

\noindent
This implies that
$\kap$ has singularities only on the curves $C_1,...,C_k$ and the singularity
of $\kap$ at the generic point of $C_j$ is of type $\lam_j$. In other words, this implies that
$$
(\calF_1,\calF_2,\kap)\in G(\calB)\backslash Y_{C_1,...,C_k}^{\lam_1,...,\lam_k}/G(\calB),
$$
i.e. we have $X_{C_1}^{\lam_1}\cdot ...\cdot X_{C_k}^{\lam_k}\subset Y_{C_1,...,C_k}^{\lam_1,...,\lam_k}$.
Thus a) follows from the fact that $X_{C_1,...,C_k}^{\lam_1,...,\lam_k}\subset Y_{C_1,...,C_k}^{\lam_1,...,\lam_k}$ and
from the fact that
all the $Y_{C_1,...,C_k}^{\lam_1,...,\lam_k}$ are disjoint.
\epr
Let us explain why \refl{cosets} implies what we need. Let us fix a finite set $\calP$ of collections
$(C_1,\lam_1),...,(C_k,\lam_k))$ with $C_1,...,C_k$ pairwise distinct and with all
$\lam_j\in\Lam_+$ and $\lam_j\neq 0$ if $k\neq 0$. Define a partial order on the set of all such collections by
declaring  that
$(C_1,\lam_1),...,(C_k,\lam_k))\geq (C_1',\lam_1',...,(C_l',\lam_l'))$ if
$k\geq l$ and after possible renumbering of the $C_j'$'s we have $C_j'=C_j$ and $\lam_j'\leq \lam_j$
for all $j=1,...,l$.We shall assume that $\calP$
is closed under the operation of replacing some $((C_1,\lam_1),...,(C_k,\lam_k))$ by a smaller element
(with respect to the above partial order). We want to show that the corresponding elements
$\alp_{C_1}(\lam_1)...\alp_{C_k}(\lam_k)\in\hsp_{crit}$ are linearly independent.

Choose $n>0$ so that the assertions (2) and (3) hold for it.
Consider now the corresponding elements
$\alp_{C_1}^{\VV_n}(\lam_1)...\alp_{C_k}^{\VV_n}(\lam_k)$. We claim that they are linearly independent as elements of
$H(G_n,G_{n,\calO})$. Indeed, let us consider the restriction map
from $H(G_n,G_{n,\calO})$ to the space $\calU_n$ of $G_{n,\calO}$-biinvariant functions on
$$
\bigcup X_{n,C_1,...,C_k}^{\lam_1,...,\lam_k}.
$$
Let also
$T_{n,C_1,...,C_k}^{\lam_1,...,\lam_k}$ be the characteristic function of
$X_{n,C_1,...,C_k}^{\lam_1,...,\lam_k}$ and set $\calU_{n,\calP}$ to be the span of
all the $T_{n,C_1,...,C_k}^{\lam_1,...,\lam_k}$ with
$(C_1,\lam_1),...,(C_k,\lam_k))\in \calP$. Then it is clear from (3) above that they form a basis of $\calU_{n,P}$.
On the other hand, it follows from (1) and (2) that  the images (in $\calU_n$) of all the
$\alp_{C_1}^{\VV_n}(\lam_1)...\alp_{C_k}^{\VV_n}(\lam_k)$ lie in $\calU_{n,\calP}$ and that
the transformation
taking $T_{n,C_1,...,C_k}^{\lam_1,...,\lam_k}$ to the image $\alp_{C_1}^{\VV_n}(\lam_1)...\alp_{C_k}^{\VV_n}(\lam_k)$
in $\calU_n$
is upper triangular with respect to the partial order on $\calP$ defined above with non-zeros on the diagonal.
This shows that the elements $\alp_{C_1}^{\VV_n}(\lam_1)...\alp_{C_k}^{\VV_n}(\lam_k)$ are linearly
independent.


\ssec{density}{Density conjecture}
The above results show that there exists an embedding
$\calR[\Lam_+]\hookrightarrow \hsp_{crit}$. Though we do not expect this map to be an isomorphism
we still would like to say that the left hand side is dense in the right
hand side in some sense. One of the ways to do this is as follows.
\conj{density}
For every $n\geq 0$
the images of $\calR[\Lam_+]$ and of $\hsp_{crit}$
in $\End(\calS_{crit}(\hGG/\Gam_n)^{G(\calA)})$ coincide.
\econj
This conjecture is obvious for $n=0$ (more precisely it follows from \refl{gago}
that $\calS_{crit}(\hGG/\Gam_0)=\calS_{crit}(\hGG/G(\calO_{\bfF}))=\CC$ and hence there is nothing to prove).
Also, it follows from \refl{homo} that \refco{density} is true for $n=1$. We do not know how
to prove \refco{density} for $n>1$.
\ssec{iw}{The Iwahori case}Let us briefly mention what happens in the Iwahori case.
In other words we want to study the algebra $H_c(\hGG,\bfI_{\calA}^{00})$ which is by definition
the algebra of endomorphisms of the functor of $\bfI_{\calA}^{00}$-invariants on the category
$\Rep_c(\hGG)$. In this case we can't describe in full detail even a dense subalgebra in
$H_c(\hGG,\bfI_{\calA}^{00})$. However, from the results \cite{K2} and \cite{GK1} one can derive
the following (details will appear elsewhere):

\smallskip
(i) For any $C\in\calC$ one can construct an embedding
$\beta_C:\overset{..}\bfH\hookrightarrow H_c(\hGG,\bfI_{\calA}^{00})$ where
$\overset{..}\bfH$ denotes the modified Cherednik's double affine Hecke algebra as in \cite{K2}.

\smallskip
(ii) There is a natural subalgebra of $H_c(\hGG,\bfI_{\calA}^{00})$ isomorphic to $\CC[\Lam]$ which
lies in the image of every $\beta_C$.

Recall that the set $\calC$ only consists of good formal curves, i.e. curves which are transversal
to the curve $t=0$. We expect that every irreducible curve in $\Spec{\calB}$ (cf. \refss{crit}) contributes
some sort of subalgebra to $H_c(\hGG,\bfI_{\calA}^{00})$. It would be extremely interesting to describe these subalgebras
explicitly.
\sec{glob}{Some remarks on $G$-bundles on $S_X^0$}

In this Section we would like to explain some constructions related to $G$-bundles
on the "surface" $S_X^0$ that were used in \refss{global}. We are going to work over an arbitrary
ground field $k$ (instead of the finite field $\FF_q$). In what follows we fix the following data:

\medskip
1) A smooth geometrically irreducible
algebraic surface $S$ over $\FF_q$;

2) A smooth projective geometrically irreducible curve $X$ over $\FF_q$;

3) A closed embedding $i:X\to S$;

4) A "retraction" $p:S\to X$ such that $p\circ i=\id$.
\footnote{This datum is not probably not necessary}

\medskip
\noindent
We denote by $\calI_X\subset \calO_S$ the sheaf of ideals corresponding to $X$.

Also we shall assume that the derived group of $G$ is simply connected.
\ssec{}{$G$-bundles on $S_X$}Let us denote by $S_X$ the formal neighbourhood of $X$ in $S$ considered as an ind-scheme.
In other words, we set
$$
S_{X,n}= \Spec_{\calO_X} (\calO_S/\calI_X^n)\quad\text{and}
\quad S_X=\underset{\to}\lim\, S_{X,n}.
$$
By a $G$-bundle on $S_{X}$ we mean a projective system of $G$-bundles on $S_{X,n}$. In other words,
to specify a $G$-bundle $\calG$ on $S_X$ we need to specify a $G$-bundle $\calG_n$ on each
$S_{X,n}$ together with  the isomorphisms
$$
\calG_n|_{S_{X,m}}\simeq \calG_m
$$
for each $n\geq m$ satisfying the standard transitivity condition.
\ssec{}{G-bundles on $S_X^0$}In \refss{global} we were talking about $G$-bundles on $S_X^0$ where the latter
was defined as the complement of $X$ in $S_X$. Let  first explain what we mean by a $G$-bundle on $S_X^0$.
Namely, we define the category of $G$-bundles on $S_X^0$ to be
category whose objects are $G$-bundles on $S_X$ and morphisms are isomorphisms of their restrictions
on $S_X^0$ (thus, by the definition, every $G$-bundle on $S_X^0$ extends to a $G$-bundle on $S_X$).

\smallskip
\noindent
{\bf Warning.}
This definition is not good for families of $G$-bundles on $S_X^0$.

We denote by $\Bun_G(S_X^0)$ the set of isomorphism classes of $G$-bundles on $S_X^0$ with respect to the above definition.
\ssec{}{Construction of the torsor $\wbg(S_X^0)$}Let us explain the construction of the torsor
$\wbg(S_X^0)$. Since in this paper we are constructing all the spaces only set-theoretically (i.e.
we do not consider families) what we actually have to do is describe a $\Pic(X)$-torsor $\wbg(S_X^0)_{\calM}$
for every $\calM\in\Bun_G(S_X^0)$.

Consider the following functor

\medskip
\centerline{Schemes over $\FF_q$ $\longrightarrow$ Sets}

\medskip
which sends every scheme $T$ to the following data:

1) A morphism $T\to X$;

2) A $G$-bundle $\calG$ on $T\underset{X}\x S_X$;

3) An isomorphism $\calG|_{T\underset{X}\x S_X^0}\simeq \pr_2^*\calM$

\noindent
where $\pr_2:T\underset{X}\x S_X^0\to S_X^0$ denotes the natural projection.
This functor is representable by an ind-scheme $\Gr_{\calM}$ endowed with
a natural proper map $\gam:\Gr_{\calM}\to X$. It is easy to see that every fiber
of $\gam$ is non-canonically isomorphic to the corresponding affine Grassmannian of $\Gr_G$ of $G$
considered in \refs{grass}. We now define the $\tbg(S_X^0)_{\calM}$ to be the set of isomorphism
classes of line bundles $\calD_{\calM}$ on $\Gr_{\calM}$ satisfying the following property:

\medskip
For every $x\in X$ the restriction of $\calD_{\calM}$ to $\gam^{-1}(x)$ is isomorphic to
the generator $\calD$ of $\Pic(\gam^{-1}(x)\simeq \Pic(\Gr_G)$.

\medskip
\noindent
Since $\calD$ is defined uniquely
up to non-canonical isomorphism and since its group of automorphisms is just the multiplicative
group (since $\Gr_G$ is ind-proper) it follows that $\calD_{\calM}$ is defined uniquely
up to tensoring with a line bundle of the form $\gam^*\calP$ where $\calP\in\Pic(X)$.
Hence $\tbg(S_X^0)$ is a homogeneous space over $\Pic(X)$. Since the map $\gam$ has sections (this is equivalent
to saying that $\calM$ extends from $S_X^0$ to $S_X$) it follows that $\tbg(S_X^0)$ is actually a
$\Pic(X)$-torsor.
\ssec{}{Proof of \refl{action-glob}}We need to construct an action of $\tgg_W$ on $\tbg(S_X^0,W)$.
To simplify the notations let us do that for $\Bun_G$ instead of $\tbg$
(the general case is treated in a very similar manner). Let us denote by $S^0_{X\backslash W}$ the pre-image
of $X\backslash W$ in $S_X^0$. Then by arguing in a way similar to the one in the proof of
\refl{gago} we see the group $\GG_W$ may be identified with the set of isomorphism classes of
$G$-bundles on $S_X^0$ trivialized at both $S_W^0$ and $S^0_{X\backslash W}$. Hence we have
$$
\Bun_G(S_X^0)=G(\calA_W)\backslash\GG_W/G(S^0_{X\backslash W})
$$
and
$$
\Bun_G(S_X^0,W)=\GG_W/G(S^0_{X\backslash W}).
$$
The right hand side has now a natural action of $\GG_W$.

\end{document}